\documentclass[preprint,12pt]{elsarticle}
\usepackage{amssymb,amsmath,amsthm}
\usepackage{float}
\usepackage{graphicx}
\usepackage{enumerate}
\usepackage{multirow}
\usepackage{amsmath,color,tikz,cancel}
\usetikzlibrary{decorations.pathreplacing}
\usepackage{hyperref}
\usepackage{url}
\usepackage{mleftright}
\usepackage{tikz}
\usepackage{tkz-euclide}
\usetikzlibrary{shapes.geometric}
\usepackage[table]{}
\usepackage{placeins}

\usepackage{tikz}
\usetikzlibrary{decorations.pathreplacing}
\usetikzlibrary{decorations.pathreplacing,decorations.markings}

\newcommand\ceil[1]{\lceil #1 \rceil}
\newtheorem{thm}{Theorem}[section]
\newtheorem{lem}[thm]{Lemma}
\newtheorem{prop}[thm]{Proposition}
\newtheorem{defn}{Definition{\rm }}

\numberwithin{equation}{section}

\usepackage{amssymb}
\usepackage{amsmath}

\journal{arXiv}

\begin{document}

\begin{frontmatter}



\title{A general upper bound on broadcast function ${{B(n)}}$ using Kn{\"o}del graph} 



\author[concordia]{Sirma Cagil Altay}
\ead{sirmaaltay@gmail.com}

\author[concordia]{Mohammadhossein Fakharan\corref{cor1}}
\ead{mohammadhossein.fakharan@mail.concordia.ca}

\author[concordia]{Hovhannes A. Harutyunyan}
\ead{haruty@cs.concordia.ca}

\address[concordia]{Department of Computer Science and Software Engineering, Concordia University, Montreal, Quebec, Canada.}

\cortext[cor1]{Corresponding author: mohammadhossein Fakharan}

\begin{abstract}
Broadcasting in a graph is the process of transmitting a message from one vertex, the originator, to all other vertices. We consider the classical model in which an informed vertex informs one of its uninformed neighbours during each time unit. A broadcast graph on ${n}$ vertices is a graph in which broadcasting can finish in ${\lceil \log_{2} n \rceil}$ time units from any originator. A minimum broadcast graph on ${n}$ vertices is a broadcast graph that has the least possible number of edges, ${B(n)}$, over all broadcast graphs on ${n}$ vertices.
This paper presents a new general upper bound on ${B(n)}$ for all odd ${n}$ except when ${n=2^{k}-1}$ for an integer $k$. In order to achieve this bound, a vertex deletion method is applied to Kn\"odel graphs on ${n+1}$ vertices. To prove our result, we characterize all informed vertices in a Kn\"odel graph after every time unit by applying the standard and reverse dimensional broadcast scheme with any cyclic shifts from an arbitrary originator. Furthermore, we completely resolve the open problem about the existence of minimum average time broadcast graphs. 
\end{abstract}



\begin{keyword}


Broadcast Function, Dimensional Broadcast Scheme, Kn\"odel Graph.

\end{keyword}

\end{frontmatter}

\section{\bf{Introduction}}

Broadcasting is a fundamental process for information dissemination in connected networks, which are commonly modeled as undirected graphs $G = (V, E)$, where $V$ represents the vertex set and $E$ represents the edge set of the graph $G$. Classical broadcasting is organized in discrete time units, beginning with an initial vertex, referred to as the originator, holding the message. During each time unit, every informed vertex (sender) transmits the message to one of its uninformed neighbors (receiver). This process continues until all vertices in the graph are informed.

The broadcasting problem has received attention in the literature due to its practical relevance and computational difficulty. Several approximation algorithms have been developed to address the challenge of computing minimum broadcast times \cite{shahin, bermond,farley1979minimum,survey1}. Despite these efforts, the problem is generally challenging. Notably, recent results have shown that the broadcasting problem is NP-complete even for some particular classes of graphs such as cactus graphs and graphs of pathwidth 2 \cite{shahin}. 

The study of the minimum broadcast graph (mbg) has a rich history, originating with Farley et al. in 1979 \cite{farley1979minimum}. In their work, they introduced the concept of minimum broadcast graphs (mbgs) and identified the first infinite family of such graphs, consisting of hypercube graphs with $n = 2^k$ vertices, where $k$ is a positive integer. Subsequently, Park and Chwa \cite{park1994recursive} proposed another family of mbgs for graphs of order $2^k$. Kn\"odel graphs are shown to be mbgs on $n = 2^k$ \cite{knodel} and $n = 2^k - 2$ vertices by Khachatrian and Haroutunian \cite{kh90} and Dinneen et al. \cite{dff91}, independently. The reader is referred to \cite{aver, bermond,  bhlp,berm, din, garg, HL23, HL99, harut5, kh91, sac} for more details and insights. In addition, we recommend surveys \cite{survey1, survey2, survey3, hrom, hrom2} for a complete overview of the topic.

Kn\"odel graphs of order $n$ have dimensions $1, \ldots, \lfloor \log n \rfloor$, similar to hypercubes of dimension $k$ with dimensions $1, \ldots, k$. A simple broadcasting algorithm for $KG(n)$, known as the \textit{dimensional broadcast scheme} \cite{bhlp}, uses a fixed dimension at each time unit. Specifically, if $D = (\alpha_1, \alpha_2, \ldots, \alpha_k)$, then in time unit $i$, every informed vertex sends the message to its $\alpha_i$-dimensional neighbor. The permutation $D$ is considered valid if all vertices receive the message after the last time unit $k$.  

Dispute Hypercubes, some permutations $D = (\alpha_1, \ldots, \alpha_k)$ in Kn\"odel  graphs are not valid dimensional broadcast schemes. In \cite{bhlp}, it was demonstrated that broadcasting in $KG(n)$ can be carried out by following specific dimension orders.

\begin{thm}\label{allvalidpermutations}
	\cite[Theorem 1]{bhlp}
	The permutation $(1, 2, \ldots, k-1, 1)$ is a valid permutation for $KG(n)$, and any cyclic shift (or the reversal of such a cyclic shift) of $(1, 2, \ldots, k-1, 1)$ is also valid, where $k = \lceil \log n \rceil$ and $n \neq 2^k$.
\end{thm}

Since identifying mbgs seems to be particularly challenging \cite{xu25}, a series of studies have introduced various techniques for constructing broadcast graphs and establishing upper bounds on broadcast function $B(n)$ for specific intervals of $n$ (see \cite{F79,garg,HL99,survey3,kh90,knodel,L94}). Farley \cite{F79} has given the first general upper bound on ${B(n)}$ using a recursive construction. 
\begin{center}
	${B(n) \leq \frac{n \lceil \log n \rceil}{2}, \text{ \ for all ${n}$}.}$
\end{center}
An even better general upper bound on ${B(n)}$ is a direct result of the number of edges of Kn\"{o}del graphs \cite{knodel}.
\begin{center}
	${B(n) \leq \frac{n \lfloor \log n \rfloor}{2} \text{ \ for all even ${n}$}.}$
\end{center}
This result was slightly improved in \cite{H08} for odd values of $n$ by using vertex addition method. In this paper, we prove that for all $n \neq 2^{k}-1$, $B(n)\leq \left \lceil \frac{n \lfloor \log n \rfloor}{2} \right \rceil$ by extending the upper bound derived from Kn\"{o}del graphs. This upper bound for an odd $n$ is obtained by applying a vertex deletion method to $KG(n+1)$ and demonstrating that the resulting graph remains a broadcast graph. The new upper bound represents an improvement of approximately $\frac{n}{2}$ over the best known general upper bound from \cite{H08}. Another advantage of our upper bound is that it is based on a direct (adhoc) construction. This explicit and straightforward construction contrasts with Farley's upper bound \cite{F79}, which is derived from a recursive construction of broadcast graphs and lacks an immediate combinatorial interpretation for arbitrary $n$.

Note that another general upper bound \cite{HL99} shows that 
\begin{equation} \label{eq1}
B(n) \leq \left( \lceil \log n \rceil - \log\left(2^{\lceil \log n \rceil} - n\right) + 1 \right) n.
\end{equation}
While \ref{eq1} is a good bound when $\log \left(2^{\lceil \log n \rceil} - n\right) =n-c$, where $c$ is constant, but in general it is as loose as $B(n)\leq \lceil \log n \rceil n$, which is twice the bound obtained by Farley. Note that the construction used to derive Equation~\ref{eq1} includes some vertices with high degree $n - 1$, while the remaining vertices have relatively small degrees. As a result, the graph is far from being regular. In contrast, in our construction, all vertices except possibly one when $\lceil \log n \rceil$ is even have degree $\lfloor \log n \rfloor$, making it a regular graph.

It is worth mentioning that there are better bounds on $B(n)$ in \cite{bermond, garg, HL23, HL99, harut5} but only for a small interval or specific values of $n$, the best among them is $B(n)\leq \frac{n}{2}(m-\frac{k+1}{2})$ for $n=2^m-2^k$. \cite{bermond, HL99, kh90}

\section{\bf{Notations}}\label{sec:pre}
The classical model of broadcasting makes the following assumptions:
\begin{itemize}
	\item there is only one originator
	\item each call involves only one informed node and one of its uninformed neighbors
	\item each call requires one time unit
	\item one node can only participate in at most one call per time unit.
\end{itemize}

A network can be modeled as a connected simple graph $G = (V, E)$, where $V$ and $E$ denote the sets of vertices and edges, respectively. 

\begin{defn}{\rm
		The {\it broadcast scheme} of a given graph $G$ from an originator vertex $v$ is a sequence of parallel calls. Each call is represented by an edge with the direction, specifying the sender and the receiver vertices.
	}
\end{defn}

\begin{defn}{\rm
		Let $G$ be a graph on $n$ vertices and $v$ be the broadcast originator in graph $G$. $b(G, v)$ defines the minimum number of time units required to broadcast from originator $v$ in graph $G$. The {\it broadcast time} of graph $G$, $b(G) = \max \{b(G, v)|v \in V(G)\}$ defines the maximum number of time units required from any originator to broadcast in graph $G$.
	}	
\end{defn}

Note that for a graph $G$ of order $n$, since the number of informed vertices can at most double during each time unit, $b(G) \geq \ceil{\log n}$. 

\begin{defn}{\rm
		If $b(G) = \ceil{\log n}$, for a graph $G$ of order $n$, then $G$ is called a {\it broadcast graph}. It means that for every vertex $v$, we have $b(G, v) = \ceil{\log n}$. A broadcast graph $G$ of order $n$ is called {\it minimum broadcast graph} ({\it mbg}) if there does not exist any other broadcast graph $G'$ of order $n$ such that $|E(G')| < |E(G)|$. The {\it broadcast function} $B(n)$ denotes the number of edges of an mbg on $n$ vertices. 
	}	
\end{defn}

\begin{defn}{\rm
		A {\it Kn{\"o}del graph} of even order $n$ and degree $\lfloor \log n \rfloor$, first introduced by W. Kn\"odel \cite{knodel}, is denoted by $KG(n) = (V, E)$, where the vertex set is $V = \{i \mid 0 \leq i \leq n - 1\}$ and the edge set is $E = \{(i, j) \mid i + j \equiv 2^s - 1 \pmod{n}, 1 \leq s \leq \lfloor \log n \rfloor, \text{ and } i, j \in V\}$. We refer to vertex $j$ as the $s$ dimensional neighbor of vertex $i$ if $i + j \equiv 2^s - 1 \pmod{n}$. 
		
	}
\end{defn}


\section{\bf{Auxiliary results}}\label{sec:alpha}
\hfill \break
Our method is to demonstrate that by removing vertex $0$ from $KG(n)$ for $n \neq 2^k$, where $k=\ceil{\log n}$, and introducing new edges, we construct another broadcast graph on $n-1$ vertices. To achieve this, based on the originator, we utilize a specific dimensional broadcast scheme in $KG(n)$, ensuring that vertex $0$ receives the message at time unit $(k-1)$ or $k$. If 0 receives at time unit $(k-1)$, then 0 is idle during the last time unit $k$, or instead of vertex $0$, one of its neighbors sends the message to the vertex that $0$ would have sent in the last time unit $k$. 

In the next lemmas, we introduce all informed vertices of $KG(n)$ after time unit $t\leq k-1$, where $k=\ceil{\log n}$, from an arbitrary originator by utilizing the standard and reverse dimensional broadcast scheme with $i$ cyclic shifts, where $0\leq i\leq k-2$. From now on, all operations are modulo $n$, the number of vertices. 


\begin{lem}\label{basestandard0}
	Let $KG(n)$ be a Kn\"odel graph of order $n < 2^k$, where $k = \lceil \log n \rceil$, and let $u$ be an arbitrary vertex. The set of informed vertices after time unit $t\leq k-1$ by using the dimensional broadcast scheme $(1,\ldots, k-1,1)$ is
	\begin{center}
		$\{2p+u|0\leq p\leq 2^{t-1}-1\} \cup \{2s-u-1|1\leq s\leq 2^{t-1}\}$.
	\end{center}
\end{lem}

\begin{proof}
	Let $A_t= \{2p+u|0\leq p\leq 2^{t-1}-1\}$ and $B_t=\{2s-u-1|1\leq s\leq 2^{t-1}\}$ for $1\leq t\leq k-1$. In the first time unit, by utilizing dimension 1, $u$ informs $1-u$. During the second time unit, these two vertices inform $3-u$ and $2+u$, respectively. We prove by induction on $t$.\\
	The base case $t = 1$ is easily verified.\\
	Assume the claim holds for some $t < k-1$. We show that it holds for $t+1$.\\
	Each vertex $v \in A_t$ has the form $v = 2p + u$ for $0 \leq p \leq 2^{t-1} - 1$. Using dimension $t+1$, such a vertex informs:
	\begin{center}
		$2^{t+1} - v - 1 = 2^{t+1} - (2p + u) - 1 = 2(2^t - p) - u - 1 = 2s - u - 1,$
	\end{center}
	where $s = 2^t - p$ ranges from $2^{t-1} + 1$ to $2^t$. These vertices, together with the vertices of $B_t$, form $B_{t+1}$.\\
	Similarly, each vertex in $B_t$ has the form $v = 2s - u - 1$ for $1 \leq s \leq 2^{t-1}$. Using dimension $t+1$, such a vertex informs:
	\begin{center}
		$2^{t+1}-2s+u= 2(2^t-s)+u=2p+u,$
	\end{center}
	where $p = 2^t - s$ ranges from $2^{t-1} + 1$ to $2^t$. These vertices, together with the vertices of $A_t$, form $A_{t+1}$.\\
	Thus, $I_{t+1} = A_{t+1} \cup B_{t+1}$, completing the induction.
\end{proof}


In the following lemma, we analyze the standard dimensional broadcast scheme with cyclic shifts. 

\begin{lem}\label{basestandardi}
	
	Let $KG(n)$ be a Kn\"odel graph of order $n < 2^k$, where $k = \lceil \log n \rceil$, and let $u$ be an arbitrary vertex. Let $1 \leq i \leq k - 2$ be an integer. Consider the dimensional broadcast scheme $(i+1, \ldots, k-1, 1, \ldots, i, i+1)$. Then the following hold:\\
	1) The set of informed vertices after time unit $t\leq k-i-1$ is
	\begin{center}
		$ \{p2^{i+1}+u|0\leq p\leq 2^{t-1}-1\}\cup \{s2^{i+1}-u-1|1\leq s\leq 2^{t-1}\} .$
	\end{center}
	2) The set of informed vertices after time unit $t\geq k-i$ is the union of the following sets:
	\begin{align*}
	A_{t} &= \left\{ p 2^{i+1} + 2j + u \;\middle|\; 0 \leq p \leq 2^{k - i - 2} - 1,\; 0 \leq j \leq 2^{t+i-k} - 1 \right\}, \\
	B_{t} &= \left\{ s 2^{i+1} + 2j - u - 1 \;\middle|\; 1 \leq s \leq 2^{k - i - 2},\; 0 \leq j \leq 2^{t+i-k} - 1 \right\}, \\
	A_{t}' &= \left\{ -p 2^{i+1} + 2j - u - 1 \;\middle|\; 0 \leq p \leq 2^{k - i - 2} - 1,\; 1 \leq j \leq 2^{t+i-k} \right\}, \\
	B_{t}' &= \left\{ -s 2^{i+1} + 2j + u \;\middle|\; 1 \leq s \leq 2^{k - i - 2},\; 1 \leq j \leq 2^{t+i-k} \right\}.
	\end{align*}
\end{lem}

\begin{proof}
	In the first time unit, by utilizing dimension $i+1$, $u$ informs $2^{i+1}-u-1$. These two vertices during the second time unit inform $2^{i+2}-u-1$ and $2^{i+1}+u$, respectively. One can see by induction that after time unit $1\leq t\leq k-i-1$, the set of informed vertices $I_t$ equals $A_t\cup B_t$, where
	\begin{center}
		$A_{t}=\{p2^{i+1}+u|0\leq p\leq 2^{t-1}-1\} \,\,\,\text{and} \,\,\, B_{t}=\{s2^{i+1}-u-1|1\leq s\leq 2^{t-1}\}$.
	\end{center}
	This holds because $A_t$ and $B_t$ contribute to each other's growth by informing new vertices at each time unit. Therefore, Part 1) is complete.\\
	By setting $t=k-i-1$, the set of informed vertices after time unit $(k-i-1)$ is 
	\begin{center}
		$I_{k-i-1}=A_{k-i-1}\cup B_{k-i-1} = \{p2^{i+1}+u|0\leq p\leq 2^{k-i-2}-1\} \cup \{s2^{i+1}-u-1|1\leq s\leq 2^{k-i-2}\}$.
	\end{center}
	In time unit $(k-i)$, the vertices of the form $p2^{i+1}+u \in A_{k-i-1}$ and $s2^{i+1}-u-1 \in B_{k-i-1}$ inform vertices of the form $-p2^{i+1}-u+1$ and $-s2^{i+1}+u+2$, respectively. Let $A_{k-i}'$ and $B_{k-i}'$ be two sets containing these newly informed vertices, respectively. By induction, one can show that for all $t \geq k - i$, the sets $A_t$ and $A_t'$ mutually contribute to each other's growth in each time unit $(t + 1)$; similarly, this holds for $B_t$ and $B_t'$. Furthermore, by the same inductive reasoning, the correctness of Part 2) follows.
\end{proof}


In the next two lemmas, we investigate similar results for the reverse dimensional broadcast schemes. We provide brief proofs of the next two lemmas, as their proofs are similar to the previous ones. 

\begin{lem}\label{basereverse0}
	Let $KG(n)$ be a Kn\"odel graph of order $n < 2^k$, where $k = \lceil \log n \rceil$, and let $u$ be an arbitrary vertex. The set of informed vertices after time unit $t\leq k-1$ by using the dimensional broadcast scheme $(k-1,\ldots, 1, k-1)$ is 
	\begin{center}
		$\{-p2^{k-t}+u|0\leq p\leq 2^{t-1}-1\} \cup \{s2^{k-t}-u-1|1\leq s\leq 2^{t-1}\}$.
	\end{center}
\end{lem}

\begin{proof}
	In the first time unit, by utilizing dimension $k-1$, $u$ informs $2^{k-1}-u-1$. These two vertices during the second dimension inform $2^{k-2}-u-1$ and $-2^{k-2}+u$, respectively. It is easy to check the correctness of the lemma by induction on $t$.
\end{proof}


In the next lemma, we show the informed vertices after any time unit $1\leq t\leq k-1$ by using the reverse dimensional broadcast scheme $(k-1,\ldots, 1)$ with $1\leq i\leq k-2$ cyclic shifts. 

\begin{lem}\label{basereversei}
	Let $KG(n)$ be a Kn\"odel graph of order $n < 2^k$, where $k = \lceil \log n \rceil$, and let $u$ be an arbitrary vertex. Let $1 \leq i \leq k - 2$ be an integer. Consider the dimensional broadcast scheme $(k-i-1, \ldots, 1, k-1, \ldots, k-i, k-i-1)$. Then the following hold:\\
	1) The set of informed vertices after time unit $t\leq k-i-1$ is
	\begin{center}
		$\{-p2^{k-i-t}+u|0\leq p\leq 2^{t-1}-1\} \cup \{s2^{k-i-t}-u-1|1\leq s\leq 2^{t-1}\}.$
	\end{center}
	2) The set of informed vertices after time unit $t\geq k-i$ is the union of the following sets:
	\begin{align*}
	A_{t} &= \left\{-j2^{2k-t-i-1} -2p + u \;\middle|\; 0 \leq j \leq 2^{t+i-k}-1,\; 0 \leq p \leq 2^{k-i-2}-1 \right\}, \\
	B_{t} &= \left\{-j2^{2k-t-i-1}+2s-u-1 \;\middle|\; 0 \leq j \leq 2^{t+i-k}-1,\; 1\leq s\leq 2^{k-i-2}\right\}, \\
	A_{t}' &= \left\{j2^{2k-t-i-1}+2p-u-1 \;\middle|\; 1 \leq j \leq 2^{t+i-k} ,\; 0 \leq p \leq 2^{k-i-2}-1 \right\}, \\
	B_{t}' &= \left\{j2^{2k-t-i-1}-2s+u \;\middle|\; 1 \leq j \leq 2^{t+i-k},\; 1 \leq s \leq 2^{k-i-2} \right\}.
	\end{align*}
\end{lem}

\begin{proof}
	In the first time unit, by utilizing dimension $k-i-1$, $u$ informs $2^{k-i-1}-u-1$. These two vertices during the second time unit inform $2^{k-i-2}-u-1$ and $-2^{k-i-2}+u$, respectively. One can see by induction that after time unit $1\leq t\leq k-i-1$, the set of informed vertices $I_t$ equals $A_t\cup B_t$, where
	\begin{center}
		$A_{t}=\{-p2^{k-i-t}+u|0\leq p\leq 2^{t-1}-1\} \,\,\,\text{and} \,\,\, B_{t}=\{s2^{k-i-t}-u-1|1\leq s\leq 2^{t-1}\}$.
	\end{center}
	This holds because $A_t$ and $B_t$ contribute to each other's growth by informing new vertices at each time unit. Therefore, Part 1 is complete.\\
	By setting $t=k-i-1$, the set of informed vertices after time unit $(k-i-1)$ is 
	\begin{center}
		$I_{k-i-1}=A_{k-i-1}\cup B_{k-i-1} = \{-2p+u|0\leq p\leq 2^{k-i-2}-1\} \cup \{2s-u-1|1\leq s\leq 2^{k-i-2}\}$.
	\end{center}
	In time unit $(k-i)$, the vertices of the form $-2p+u \in A_{k-i-1}$ and $2s-u-1 \in B_{k-i-1}$ inform vertices of the form $2^{k-1}+2p-u-1$ and $2^{k-1}-2s+u$, respectively. Let $A_{k-i}'$ and $B_{k-i}'$ be two sets containing these newly informed vertices, respectively. By induction on $t$, one can show that for all $t \geq k - i$, the sets $A_t$ and $A_t'$ mutually contribute to each other's growth in each time unit $(t + 1)$; similarly, this holds for $B_t$ and $B_t'$. Furthermore, by the same inductive reasoning, the correctness of Part 2 follows.
\end{proof}


In the next proposition, we consider the case where the originator is an $i$-dimensional neighbor of vertex $0$, with $1 \leq i \leq k - 1$. We prove that by utilizing $i$ cyclic shifts of the reverse dimensional broadcast scheme, vertex $0$ remains uninformed after time unit $(k - 2)$.

\begin{prop}\label{propodd}
	Let $KG(n)$ be a Kn\"odel graph of order $n$ with $k=\ceil {\log n}$ and $0,u\in V(KG(n))$ be two $i$ dimensional neighbors for an integer $1\leq i\leq k-1$. If the originator and the dimensional broadcast scheme are $u$ and $(i+1,i+2, \ldots, k-1,1,\ldots,i,i+1)$, respectively, then vertex $0$ is uninformed after time unit $k-2$ and will be informed by $u$ in time unit $k-1$.  (In case of $i=k-1$, the dimensional broadcast scheme is $(1,\ldots, k-1,1)$).
\end{prop}

\begin{proof}
	First, suppose that $i\neq k-1$. In Lemma \ref{basestandardi}, let $u=2^i-1$ and $t=k-2$. If $i=1$, then by Part 1, the set of informed vertices $I$ after time unit $(k-2)$ is 
	\begin{center}
		$\{4p+1|0\leq p\leq 2^{k-3}-1\}\cup \{4s-2|1\leq s\leq 2^{k-3}\} .$
	\end{center}
	If $i\geq 2$, then by Part 2 of Lemma \ref{basestandardi}, the set of informed vertices $I$ after time unit $(k-2)$ is the union of the following sets:\\
	$A=\{(2p+1)2^i+j| -1\leq j \leq 2^{i-1}-2\}\\
	B=\{-(2p+1)2^{i}+j \text{ (mod } n)|1\leq j\leq 2^{i-1}\},\\
	$
	where $0\leq p\leq 2^{k-i-2}-1$.\\
	Now, we check if $0 \overset{n}{\equiv} n \in I$.\\
	Obviously, $0\notin I$ when $i=1$.\\
	Let $i\geq 2$. Note that $0<\min A= 2^i-1<\max A=2^{k-1}-2^{i-1}-2<n$. Hence, $0\notin A$.\\
	Since $i\neq k-1$, then $-n<\min B= -2^{k-1}+2^i+1<\max B= -2^{i-1}<0$. Thus, all members of $B$ (mod $n$) are between 0 and $n$ and so $0\notin B$. Therefore, $0\notin I$ and the proof is complete. \\
	Now, suppose that $i=k-1$. In this case, $u=2^{k-1}-1$ and the dimensional broadcast scheme is $(1,\ldots, k-1,1)$. By Lemma \ref{basestandard0}, the set of informed vertices $I$ after time unit $(k-2)$ is the union of the following sets:\\
	$A=\{2^{k-1}+2p-1| -1\leq p \leq 2^{k-3}-1\}\\
	B=\{2p-2^{k-1} \text{ (mod } n)|1\leq p\leq 2^{k-3}\}.\\
	$
	Since all members of $A$ are odd vertices, $0\notin A$. \\
	It is worth mentioning that $-n<\min B= -2^{k-1}+2<\max B=-2^{k-2}<0$. Thus, all members of $B$ modulo $n$ are between 0 and $n$. Therefore, $0\notin B$ and so $0\notin I$.\\
	Finally, it is obvious that $u$ informs 0 in time unit $k-1$ by using dimension $i$. 
\end{proof}


In the next proposition, we consider the case where the originator is an $i$-dimensional neighbor of vertex $0$, with $1 \leq i \leq k - 1$. We prove that by utilizing $k-i$ cyclic shifts of the reverse dimensional broadcast scheme, vertex $0$ remains uninformed after time unit $(k - 2)$.

\begin{prop}\label{propeven}
	Let $KG(n)$ be a Kn\"odel graph of order $n$ with $k=\ceil {\log n}$ and $0,u\in V(KG(n))$ be two $i$ dimensional neighbors for an integer $1\leq i\leq k-1$. If the originator and the dimensional broadcast scheme are $u$ and $(i-1,i-2, \ldots,1,k-1, k-2,\ldots,i,i-1)$, respectively, then vertex $0$ is uninformed after time unit $k-2$ and will be informed by $u$ in time unit $k-1$.  (In case of $i=1$, the dimensional broadcast scheme is $(k-1,\ldots, 1,k-1)$).
\end{prop}

\begin{proof}
	First, suppose that $i\neq 1$. In Lemma \ref{basereversei}, consider $k-i$ cyclic shifts instead of $i$ and let $u=2^i-1$ and $t=k-2$.\\
	If $i=k-1$, then by Part 1, the set of informed vertices $I$ after time unit $(k-2)$ is 
	\begin{center}
		$\{2^{k-1}-2p-1|0\leq p\leq 2^{k-3}-1\} \cup \{2s-2^{k-1}|1\leq s\leq 2^{k-3}\}.$
	\end{center} 
	If $i\leq k-2$, then by Part 2 of Lemma \ref{basereversei}, the set of informed vertices $I$ after time unit $(k-2)$ is the union of the following sets:\\
	\begin{align*}
	A_{k-2} &= \left\{-j2^{i+1} -2p + 2^i-1 \;\middle|\; 0 \leq j \leq 2^{k-i-2}-1,\; 0 \leq p \leq 2^{i-2}-1 \right\}, \\
	B_{k-2} &= \left\{-j2^{i+1}+2s-2^i \;\middle|\; 0 \leq j \leq 2^{k-i-2}-1,\; 1\leq s\leq 2^{i-2}\right\}, \\
	A_{k-2}' &= \left\{j2^{i+1}+2p-2^i \;\middle|\; 1 \leq j \leq 2^{k-i-2} ,\; 0 \leq p \leq 2^{i-2}-1 \right\}, \\
	B_{k-2}' &= \left\{j2^{i+1}-2s+2^i-1 \;\middle|\; 1 \leq j \leq 2^{k-i-2},\; 1 \leq s \leq 2^{i-2} \right\}.
	\end{align*}
	Now, we check if $0 \overset{n}{\equiv} n \in I$.\\
	Obviously, $0\notin I$ when $i=k-1$.\\
	Let $i\leq k-2$. Note that $-n<\min B_{k-2}= -2^{k-1}+2^i+2<\max B_{k-2}=-2^{i-1}<0$. Hence, all members of $B_{k-2}$ (mod $n$) are between 0 and $n$ and so $0\notin B_{k-2}$.\\
	Also, $0<\min A_{k-2}'= 2^i<\max A_{k-2}'= 2^{k-1}-2^{i-1}-2<n$. Thus,  $0\notin A_{k-2}'$. Therefore, $0\notin I$ and the proof is complete. \\
	Now, suppose that $i=1$. In this case, $u=1$ and the dimensional broadcast scheme is $(k-1,\ldots,1, k-1)$. By Lemma \ref{basereverse0}, the set of informed vertices $I$ after time unit $(k-2)$ is 
	\begin{center}
		$\{-4p+1|0\leq p\leq 2^{k-3}-1\} \cup \{4s-2|1\leq s\leq 2^{k-3}\}$.
	\end{center}
	It is easy to check that $0\notin I$. \\
	Finally, it is obvious that $u$ informs 0 in time unit $k-1$ by using dimension $i$.
\end{proof}


In the next claim, we consider the case where the originator $u$ is a vertex that is not a neighbor of vertex 0. We prove that based on $u$, there exists a dimensional broadcast scheme such that the vertex 0 is either uninformed after time unit $(k-1)$ or remains uninformed after time unit $(k-2)$, receives the message in time unit $(k-1)$, and remains idle during the last time unit $k$.

\begin{prop}\label{propother}
	Let $KG(n)$ be a Kn\"odel graph of order $n$ with $k=\ceil {\log n}$. If $u$ is an arbitrary vertex that is not adjacent to vertex 0, then there exists a dimensional broadcast scheme such that vertex 0 is either uninformed after time unit $(k-1)$ or remains uninformed after time unit $(k-2)$, receives the message in time unit $(k-1)$, and remains idle during the last time unit $k$.
\end{prop}

\begin{proof}
	There are two cases:\\
	{\bf{Case 1.} } Assume that $u$ is odd. In this case, we show that there exists a dimensional broadcast scheme, dependent on $u$, such that vertex 0 remains uninformed after time unit $(k-1)$.\\
	If $u\geq 2^{k-1}$, then consider dimensional broadcast scheme $(1,\ldots,k-1,1)$.  By Lemma \ref{basestandard0}, the set of informed vertices $I_{k-1}$ after time unit $(k-1)$ equals to the union of the following sets: 
	\begin{center}
		$A_{k-1} = \{2p+u|0\leq p\leq 2^{k-2}-1\} \,\,\, \text{and} \,\,\, B_{k-1} = \{2s-u-1|1\leq s\leq 2^{k-2}\}$.
	\end{center}
	Since $u$ is odd, $0\notin A_{k-1}$. Also, $-n<\min B_{k-1}= 1-u<\max B_{k-1}= 2^{k-1}-u-1<0$. Hence, $0\notin B_{k-1}$ and so $0\notin I_{k-1}$. \\
	If $2^i \leq u \leq 2^{i+1} - 3$ for some integer $2 \leq i \leq k - 2$, then consider the broadcast scheme obtained by applying $i$ cyclic shifts to the standard dimensional broadcast order, i.e., $(i + 1, \ldots, k - 1, 1, \ldots, i, i + 1)$.
	By Lemma \ref{basestandardi}, the set of informed vertices $I_{k-1}$ after time unit $(k-1)$ is the union of the following sets:
	\begin{align*}
	A_{k-1} &= \left\{ p 2^{i+1} + 2j + u \;\middle|\; 0 \leq p \leq 2^{k - i - 2} - 1,\; 0 \leq j \leq 2^{i-1} - 1 \right\}, \\
	B_{k-1} &= \left\{ s 2^{i+1} + 2j - u - 1 \;\middle|\; 1 \leq s \leq 2^{k - i - 2},\; 0 \leq j \leq 2^{i-1} - 1 \right\}, \\
	A_{k-1}' &= \left\{ -p 2^{i+1} + 2j - u - 1 \;\middle|\; 0 \leq p \leq 2^{k - i - 2} - 1,\; 1 \leq j \leq 2^{i-1} \right\}, \\
	B_{k-1}' &= \left\{ -s 2^{i+1} + 2j + u \;\middle|\; 1 \leq s \leq 2^{k - i - 2},\; 1 \leq j \leq 2^{i-1} \right\}.
	\end{align*}
	Since $u$ is odd, $0\notin A_{k-1}, B_{k-1}'$.\\
	Suppose that $0\in B_{k-1}$. Thus, there exist $1 \leq s \leq 2^{k - i - 2}$ and $0 \leq j \leq 2^{i-1} - 1$ such that $s 2^{i+1} + 2j - u - 1\overset{n}{=}0$. If $s 2^{i+1} + 2j - u - 1=0$, then $u\geq 2^{i+1}-1$, a contradiction. If $s 2^{i+1} + 2j - u - 1 =n$, then $u\leq 2^{k-1} +2^i-n-3< 2^i$, a contradiction. Hence, $0\notin B_{k-1}$. \\ 
	Assume that $0\in A_{k-1}'$. Then there exist $0 \leq p \leq 2^{k - i - 2} - 1$ and $1 \leq j \leq 2^{i-1}$ such that $-p 2^{i+1} + 2j - u - 1\overset{n}{=}0$. If $-p 2^{i+1} + 2j - u - 1=0$, then $u\leq 2^{i}-1$, a contradiction. If $-p 2^{i+1} + 2j - u - 1 =-n$, then $u\geq n-2^{k-1} +2^{i+1}-1> 2^{i+1}$, a contradiction. Hence, $0\notin A_{k-1}'$. \\ 
	Thus, $0\notin I_{k-1}$. Note that since $u$ is not adjacent to vertex 0, then $u\neq 2^i-1$ for $1\leq i\leq k-1$. Therefore, for every odd vertex $u$ (which is not adjacent to 0), there exists a dimensional broadcast scheme such that 0 remains uninformed after time unit $(k-1)$.
	\vspace{3mm}
	
	\noindent {\bf{Case 2.} } Assume that $u$ is even. There are 5 sub-cases.\\
	Sub-case 1. $u\leq n-2^{k-1}$. Consider the dimensional broadcast scheme $(1,\ldots,k-1,1)$.  By Lemma \ref{basestandard0}, the set of informed vertices $I_{k-1}$ after time unit $(k-1)$ equals to the union of the following sets: 
	\begin{center}
		$A_{k-1} = \{2p+u|0\leq p\leq 2^{k-2}-1\} \,\,\, \text{and} \,\,\, B_{k-1} = \{2s-u-1|1\leq s\leq 2^{k-2}\}$.
	\end{center}
	Since $u$ is even, $0\notin B_{k-1}$. Also, $0<\min A_{k-1}= u<\max A_{k-1}= 2^{k-1}+u-2<n$. Hence, $0\notin A_{k-1}$ and so $0\notin I_{k-1}$. 
	\vspace{3mm}
	
	\noindent Sub-case 2. $u< 2^{k-2}$. Consider the broadcast scheme obtained by applying $k-2$ cyclic shifts to the standard dimensional broadcast scheme, i.e. $(k-1,1,\ldots,k-2,k-1)$.  By setting $i=k-2$ in Lemma \ref{basestandardi}, the set of informed vertices $I_{k-1}$ after time unit $(k-1)$ equals to the union of the following sets: 
	\begin{align*}
	A_{k-1} &= \left\{ 2j + u \;\middle|\;  0 \leq j \leq 2^{k-3} - 1 \right\}, \\
	B_{k-1} &= \left\{ 2^{k-1} + 2j - u - 1 \;\middle|\; 0 \leq j \leq 2^{k-3} - 1 \right\}, \\
	A_{k-1}' &= \left\{ 2j - u - 1 \;\middle|\;1 \leq j \leq 2^{k-3} \right\}, \\
	B_{k-1}' &= \left\{ -2^{k-1} + 2j + u \;\middle|\; 1 \leq j \leq 2^{k-3} \right\}.
	\end{align*}
	Since $u$ is even, $0\notin B_{k-1}, A_{k-1}'$.\\
	Suppose that $0\in A_{k-1}$. Thus, there exists $0 \leq j \leq 2^{k-3} - 1$ such that $u+2j\overset{n}{=}0$. So, $u+2j=n$. Hence, $u\geq n-2^{k-2}+2>2^{k-2}$, a contradiction. Therefore, $0\notin A_{k-1}$. \\ 
	Assume that $0\in B_{k-1}'$. Then there exists $1 \leq j \leq 2^{k-3}$ such that $-2^{k-1} + 2j + u\overset{n}{=}0$. So, $-2^{k-1} + 2j + u=0$. Hence, $u\geq 2^{k-2}$, a contradiction. Thus, $0\notin B_{k-1}'$.\\
	Therefore, $0\notin I_{k-1}$. 
	\vspace{3mm}
	
	\noindent Sub-case 3. $u\geq 2^{k-1}$. Consider the dimensional broadcast scheme $(k-1,\ldots,1,k-1)$.  By Lemma \ref{basereverse0}, the set of informed vertices $I_{k-1}$ after time unit $(k-1)$ equals to the union of the following sets: 
	\begin{center}
		$A_{k-1} = \{-2p+u|0\leq p\leq 2^{k-2}-1\} \,\,\, \text{and} \,\,\, B_{k-1} =\{2s-u-1|1\leq s\leq 2^{k-2}\}$.
	\end{center}
	Since $u$ is even, $0\notin B_{k-1}$. Also, $0<\min A_{k-1}=-2^{k-1}+u+2 <\max A_{k-1}= u<n$. Hence, $0\notin A_{k-1}$ and so $0\notin I_{k-1}$. 
	\vspace{3mm}
	
	\noindent Sub-case 4. $u> n-2^{k-2}$. Consider the broadcast scheme obtained by applying $1$ cyclic shift to the reverse dimensional broadcast scheme, i.e. $(k-2,\ldots,1,k-1,k-2)$.  By setting $i=1$ in Lemma \ref{basereversei}, the set of informed vertices $I_{k-1}$ after time unit $(k-1)$ equals to the union of the following sets: 
	\begin{align*}
	A_{k-1} &= \left\{ -2p + u \;\middle|\;  0 \leq p \leq 2^{k-3} - 1 \right\}, \\
	B_{k-1} &= \left\{ 2s - u - 1 \;\middle|\; 1 \leq s \leq 2^{k-3} \right\}, \\
	A_{k-1}' &= \left\{ 2^{k-1} +2p - u - 1\;\middle|\;0\leq p\leq 2^{k-3}-1 \right\},\\
	B_{k-1}' &= \left\{ 2^{k-1} -2s + u \;\middle|\; 1 \leq s \leq 2^{k-3} \right\}.
	\end{align*}
	Since $u$ is even, $0\notin B_{k-1}, A_{k-1}'$.\\
	Suppose that $0\in A_{k-1}$. Thus, there exists $0 \leq p \leq 2^{k-3} - 1$ such that $u-2p\overset{n}{=}0$. So, $u-2p=0$. Hence, $u\leq 2^{k-2}-2<n-2^{k-2}$, a contradiction. Therefore, $0\notin A_{k-1}$. \\ 
	Assume that $0\in B_{k-1}'$. Then there exists $1 \leq s \leq 2^{k-3}$ such that $2^{k-1} -2s + u\overset{n}{=}0$. So, $2^{k-1} -2s + u=n$. Hence, $u\leq n-2^{k-2}$, a contradiction. Thus, $0\notin B_{k-1}'$.\\
	Therefore, $0\notin I_{k-1}$. \\
	Note that in all previous subcases, there exists a dimensional broadcast scheme, depending on $u$, under which vertex $0$ remains uninformed after time unit $(k - 1)$. However, for the next subcase, we provide a dimensional broadcast scheme such that 0 remains uninformed after time unit $(k-2)$, receives the message during time unit $(k-1)$, and remains idle during the last time unit $k$.
	\vspace{3mm}
	
	\noindent Sub-case 5. $n-2^{k-1}+2\leq u\leq \min \{2^{k-1}-2,n-2^{k-2}\}$. In this case, consider the dimensional broadcast scheme $(1,\ldots,k-1,1)$.  By Lemma \ref{basestandard0}, the set of informed vertices $I_{k-2}$ after time unit $(k-2)$ equals to the union of the following sets: 
	\begin{center}
		$A_{k-2} = \{2p+u|0\leq p\leq 2^{k-3}-1\} \,\,\, \text{and} \,\,\, B_{k-2} = \{2s-u-1|1\leq s\leq 2^{k-3}\}$.
	\end{center}
	Since $u$ is even, $0\notin B_{k-2}$. Also, $0<\min A_{k-2}= u<\max A_{k-2}= 2^{k-2}+u-2<n$. Hence, $0\notin A_{k-2}$ and so $0\notin I_{k-2}$. It is worth mentioning that for $s=\dfrac{u-n+2^{k-1}}{2}$, $2s-u-1 \overset{n}{=} 2^{k-1}-1$. Thus, $2^{k-1}-1\in I_{k-1}$, and so $0\in I_{k-1}$. Now, we claim that $1\in I_{k-1}$. Let $s=\frac{u+2}{2}$. We have $2s-u-1=1$ and $s\leq 2^{k-2}$. Hence, $1\in I_{k-1}$. Thus, 0 remains idle during the last time unit.  \\
	Hence, for every even originator vertex, there exists a dimensional broadcast scheme such that 0 remains idle after time unit $(k-2)$ or remains idle after time unit $(k-2)$, receives the message during time unit $(k-1)$, and remains idle during the last time unit. Therefore, the proof is complete. 
\end{proof}

\section{Construction of a Broadcast Graph}
In this section, we remove one vertex, say vertex $0$, from $KG(n)$ and add $\ceil{\frac{k-1}{2}}$ edges, where $k=\ceil{\log n}$, to construct a new broadcast graph on $n - 1$ vertices. Using this method, we construct a broadcast graph for every odd number of vertices. Consequently, we derive an upper bound for $B(n)$ for all integers $n$.

In general, when a vertex $v$ is deleted from a broadcast graph, one must add $\binom{\deg(v)}{2}$ edges (forming a clique among $v$'s neighbors) to ensure that the resulting graph on $n - 1$ vertices remains a broadcast graph. However, in $KG(n)$, it suffices to add only $\left\lceil \frac{\deg(v)}{2} \right\rceil$ edges, thanks to its good structural properties.

\begin{figure}[H]
	\begin{center}
	\resizebox{10cm}{!}{
		\begin{tikzpicture}[node distance =2 cm and 2cm, thin, state/.style ={circle, fill=black, draw, black, minimum width =0.3 cm}]
		
		\tikzset{edge/.style = {->,> = latex'}}
		
		\node[draw, circle, scale=0.5] (10) at (-2,.5) {}; 
		\node (u) at (-2,.7)[scale =.4] {$0$};
		
		\node[draw, circle, scale=0.5] (a1) at (-4.8,-1) {};
		\node (u) at (-4.8,-1.3)[scale =.4] {$1$};
		
		\node[draw, circle, scale=0.5] (a2) at (-3.8,-1) {};
		\node (u) at (-3.8,-1.3)[scale =.4] {$3$};
		
		\node[draw, circle, scale=0.5] (a3) at (-2.8,-1) {};
		\node (u) at (-2.8,-1.3)[scale =.4] {$7$};
		
		\node[draw, circle, scale=0.5] (a4) at (-1.8,-1) {};
		\node (u) at (-1.8,-1.3)[scale =.4] {$15$};
		

		\node[draw, circle, scale=0.5] (a5) at (-.8,-1) {};
		\node (u) at (-0.8,-1.3)[scale =.4] {$2^{k-3}-1$};

		\node[draw, circle, scale=0.5] (a6) at (.2,-1) {};
		\node (u) at (0.2,-1.3)[scale =.4] {$2^{k-2}-1$};
		
		\node[draw, circle, scale=0.5] (a7) at (1.2,-1) {};
		\node (u) at (1.2,-1.3)[scale =.4] {$2^{k-1}-1$};
		\draw[dashed] (10) -- (a1);
		\draw[dashed] (10) -- (a2);
		\draw[dashed] (10) -- (a3);
		\draw[dashed] (10) -- (a4);
		\draw[dashed] (10) -- (a5);
		\draw[dashed] (10) -- (a6);
		\draw[dashed] (10) -- (a7);
		\draw (a1) -- (a2);
		\draw (a3) -- (a4);
		\draw (a5) -- (a6);
		\draw (a6) -- (a7);
		\node at (-1.3,-1) {\small $\cdots$};
		\end{tikzpicture}}

	\caption{Subgraph of $KG'(n-1)$ when $k-1$ is odd. Dashed lines represent removed edges, and solid lines represent newly added edges.}
\end{center}
	\label{construction2}
\end{figure}
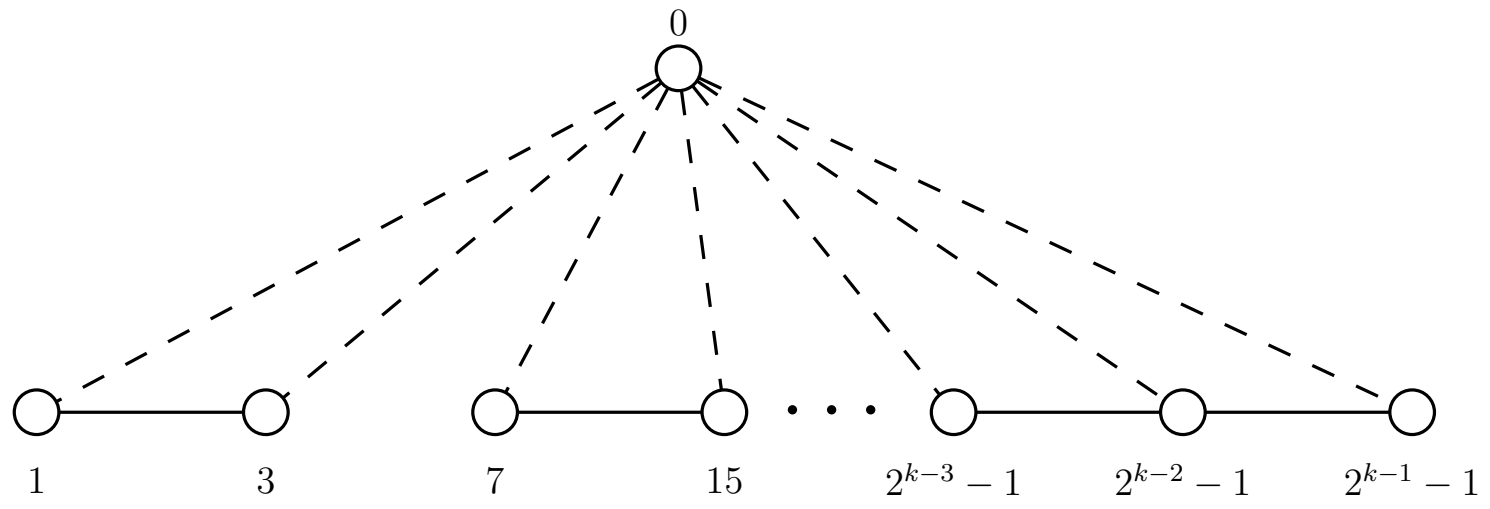

Let $KG(n)$ be a Kn\"odel graph of order $n$ with $k=\ceil {\log n}$ and $n\neq 2^k$. We introduce a new graph $KG'(n-1)$ by removing vertex 0 and all its incident edges, and adding edges $(2^{i-1}-1,2^i-1)$ for all even $2\leq i\leq k-1$. If $k-1$ is odd, then we also add edge $(2^{k-2}-1,2^{k-1}-1)$ as shown in Figures \ref{construction1} and \ref{construction2}. 

\begin{figure}[H]
	\begin{center}

	\resizebox{10cm}{!}{
		\begin{tikzpicture}[node distance =2 cm and 2cm, thin, state/.style ={circle, fill=black, draw, black, minimum width =0.3 cm}]
		
		\tikzset{edge/.style = {->,> = latex'}}
		
		\node[draw, circle, scale=0.5] (10) at (-2,.5)[minimum size=0.02cm] {};
		\node (u) at (-2,.7)[scale =.4] {$0$};
		\node[draw, circle, scale=0.5] (a1) at (-4.5,-1)[minimum size=0.01cm] {};
		\node (u) at (-4.5,-1.3)[scale =.4] {$1$};
		\node[draw, circle, scale=0.5] (a2) at (-3.5,-1)[minimum size=0.01cm] {};
		\node (u) at (-3.5,-1.3)[scale =.4] {$3$};
		\node[draw, circle, scale=0.5] (a3) at (-2.5,-1)[minimum size=0.01cm] {};
		\node (u) at (-2.5,-1.3)[scale =.4] {$7$};
		\node[draw, circle, scale=0.5] (a4) at (-1.5,-1)[minimum size=0.01cm] {};
		\node (u) at (-1.5,-1.3)[scale =.4] {$15$};
		\node[draw, circle, scale=0.5] (a5) at (-.5,-1)[minimum size=0.01cm] {};
		\node (u) at (-0.5,-1.3)[scale =.4] {$2^{k-2}-1$};
		\node[draw, circle, scale=0.5] (a6) at (.5,-1)[minimum size=0.01cm] {};
		\node (u) at (0.5,-1.3)[scale =.4] {$2^{k-1}-1$};
		\draw[dashed] (10) -- (a1);
		\draw[dashed] (10) -- (a2);
		\draw[dashed] (10) -- (a3);
		\draw[dashed] (10) -- (a4);
		\draw[dashed] (10) -- (a5);
		\draw[dashed] (10) -- (a6);
		\draw (a1) -- (a2);
		\draw (a3) -- (a4);
		\draw (a5) -- (a6);
		\node at (-1,-1) {\small $\cdots$};
		\end{tikzpicture}}
	\end{center}
	\caption{Subgraph of $KG'(n-1)$ when $k-1$ is even. Dashed lines represent removed edges, and solid lines represent newly added edges.}
	\label{construction1}
\end{figure}
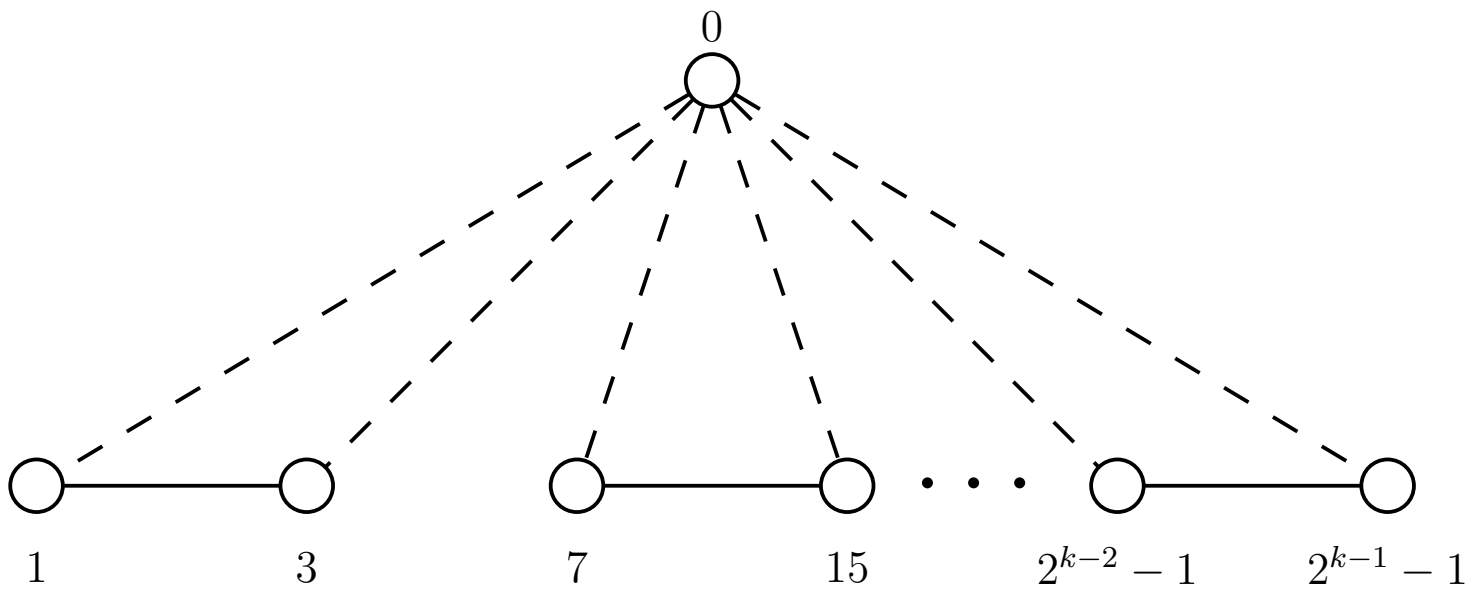 

In the next theorem, we prove that $KG'(n-1)$ is a broadcast graph. 

\begin{thm}\label{broad}
	Let $KG(n)$ be a Kn\"odel graph of order $n$ with $k=\ceil {\log n}$ such that $n\neq 2^k$. The constructed graph $KG'(n-1)$ from $KG(n)$ is a broadcast graph.  
\end{thm}

\begin{proof}
	If the originator is $2^i - 1$ for some odd integer $1 \leq i \leq k - 2$, we claim that the dimensional broadcast scheme $(i+1, \ldots, k-1, 1, \ldots, i, i+1)$ is valid for graph $KG'(n-1)$. This holds because, in $KG(n)$, by Proposition \ref{propodd}, vertex 0 remains uninformed after time unit $(k - 2)$, receives the message from the originator at time unit $(k - 1)$, and then informs vertex $2^{i+1} - 1$ in the last time unit. In the constructed graph $KG'(n - 1)$, the originator directly informs vertex $2^{i+1} - 1$ at time unit $(k - 1)$, replacing the role of vertex 0 in graph $KG(n)$. The rest of the broadcast scheme remains unchanged. Therefore, every vertex will be informed after the last time unit $k$.\\
	If the originator is $2^i - 1$ for some even integer $1 \leq i \leq k - 1$, a similar argument applies using the dimensional broadcast scheme $(i-1, \ldots, 1, k-1, \ldots, i, i-1)$. In $KG(n)$, by Proposition \ref{propeven}, vertex 0 receives the message at time unit $(k - 1)$ and informs vertex $2^{i-1} - 1$ at time $k$. In $KG'(n - 1)$, the originator directly informs $2^{i-1} - 1$ at time unit $(k - 1)$, ensuring all vertices are informed after the last time unit $k$.\\
	If $k-1$ is odd, then similar to the previous case, the dimensional broadcast scheme $(k-2, \ldots, 1, k-1, k-2)$ is valid for $KG'(n-1)$. \\
	If the originator is a vertex that is not adjacent to 0, then, in $KG(n)$, by Proposition \ref{propother}, there exists a dimensional broadcast scheme $D$ such that 0 remains uninformed after time unit $(k-1)$ or remains uninformed after time unit $(k-2)$, receives the message in time unit $(k-1)$ and remains idle during the last time unit. Since vertex 0 does not contribute to the broadcast by using $D$, removing it does not affect the dissemination of the message in $KG'(n-1)$. Therefore, $KG'(n-1)$ is a broadcast graph. 
\end{proof}

\begin{thm}\label{bound}
	For all $n \neq 2^{k}-1$ where $k \geq 2$ 
	\begin{center}
		$B(n) \leq \left \lceil \dfrac{n \lfloor \log n \rfloor}{2} \right \rceil$.  
	\end{center}
\end{thm}

\begin{proof}
	If $n$ is even, then the number of edges in $KG(n)$ is $\dfrac{n \lfloor \log n \rfloor}{2}$, and there is nothing to prove. Suppose that $n$ is odd. Consider graph $KG'(n)$ obtained from $KG(n+1)$. The number of edges in $KG'(n)$ is $\dfrac{n(k - 1)}{2}$ if $k - 1$ is even, and $\dfrac{n(k - 1) + 1}{2}$ if $k - 1$ is odd. Since by Theorem \ref{broad}, $KG'(n)$ is a broadcast graph, then 
	\begin{center}
		$B(n) \leq \left \lceil \dfrac{n \lfloor \log n \rfloor}{2} \right \rceil$. 
	\end{center}
\end{proof}

In the case $n=2^{k}-1$, it is proved in \cite{L94} that for all ${k \geq 2}$
$${B(2^{k}-1) \geq \frac{(k-1)(2^{k}-1)}{2} + \frac{2^{k}-1}{2(k+1)}.}$$
This lower bound is greater than $\left \lceil \dfrac{n \lfloor \log n \rfloor}{2} \right \rceil$.
Hence the general upper bound ${\left \lceil \dfrac{n \lfloor \log n \rfloor}{2} \right \rceil}$ cannot be expanded to the case ${n=2^{k}-1}$.

\section{Minimum Average Broadcast Graphs}

The concept of {\it minimum average broadcast graphs} has been investigated in \cite{liest}. Its definition is as follows: A broadcast scheme can be represented by a spanning tree $T_u$ of a graph $G$, rooted at the originator $u$, where each vertex $v$ is assigned a value $t(v)$ indicating the time at which $v$ receives the message. We assume $t(u) = 0$. The \emph{average broadcast time} of a scheme $T_u$ is defined as
\begin{center}
	$ b_a(T_u) = \dfrac{1}{|V|} \displaystyle\sum_{v \in V} t(v). $
\end{center}
For a fixed originator $u$, the \emph{minimum average broadcast time for $u$} is defined as
\begin{center}
	$ b_a(u) = \min \{ b_a(T_u) \}, $
\end{center}
where the minimum is taken over all broadcast trees rooted at $u$. The \emph{minimum average broadcast time for a graph} $G$ is then defined as
\begin{center}
	$ b_a(G) = \displaystyle\max_{u \in V(G)} \{ b_a(u) \}. $
\end{center}
Since the complete graph $K_n$ has the minimum $b_a(G)$ among all graphs on $n$ vertices, we define $\beta_a(n) = b_a(K_n)$, and refer to it as the \emph{minimum average broadcast time for $n$ vertices}.

A \emph{minimum average time broadcast graph} (matbg) is a graph $G$ on $n$ vertices such that $b_a(G) = \beta_a(n)$, and for any other graph $G'$ on $n$ vertices with $b_a(G') = \beta_a(n)$, the graph $G$ has no more edges than $G'$. In other words, matbg's are graphs on $n$ vertices with the fewest edges that still allow broadcasting in minimum average time from any originator. Note that in order to achieve the minimum average broadcast time $\beta_a(n)$, the number of informed vertices must double in each time unit except the last time unit. Consequently, the degree of every vertex $v$ in a minimum average time broadcast graph (matbg) on $n$ vertices must satisfy $\deg(v) \geq \lfloor \log n \rfloor$.
We denote by $B_a(n)$ the number of edges in a matbg on $n$ vertices.

note that Kn\"odel graphs $KG(n)$ are matbgs for all even values of $n$. In Theorems 3.9 and 3.11 of \cite{liest}, it was shown that for $2^k + 1 \leq n \leq 2^k + 2^{k-1}$ with $k \geq 5$, we have $B_a(n) = \left\lceil \frac{kn}{2} \right\rceil$, by explicitly constructing a matbg on $n$ vertices using the multiple fixed step graph method. This construction cannot be extended for the second half of the interval $(2^{k-1},2^{k})$.
While the authors successfully addressed the case for the first half of the interval $(2^{k-1}, 2^k)$, the broader problem of finding matbg graphs throughout the entire interval remained an open problem. Theorem \ref{bound} provides a new construction that also yields a matbg. As a result, we completely resolve the open problem by extending the characterization of matbg graphs to the full interval.

\FloatBarrier







\section{Conclusion}
In this paper, we established a new general upper bound on the broadcast function $B(n)$ for all odd integers $n$, except when $n = 2^k - 1$ for some integer $k$. Our approach is based on a vertex deletion technique applied to Kn\"odel graphs on $n+1$ vertices.

To support our main result, we developed a detailed characterization of the informed vertices in Kn\"odel graphs at each time unit, using both the standard and reverse dimensional broadcast schemes with arbitrary cyclic shifts from an arbitrary originator. This analysis provides a comprehensive view of the broadcast schemes within these graphs and helps us better understand their structural properties.

Moreover, we addressed and completely resolved an open problem regarding the existence of minimum average time broadcast graphs, thereby contributing to the broader theory of information dissemination in networks. Our findings not only improve known bounds but also offer new tools and perspectives for future investigations.

\end{document}